# A polygonal scaled boundary finite element method for solving heat conduction problems

Yang Yang [1,2,*], Zongliang Zhang [1,*], Yelin Feng [1], Yuzhen Yu [2], Kun Wang [1,2], and Lihui Liang [1]

[1] PowerChina Kunming Engineering Corporation Limited, Kunming 650051, China;

[2] Department of Hydraulic Engineering, Tsinghua University, Beijing 100084, China;

\* Correspondence: Yang Yang, yangyhhu@foxmail.com; Zongliang Zhang, smiba@qq.com;

**Abstract**

This paper presents a steady-state and transient heat conduction analysis framework using the polygonal scaled boundary finite element method (PSBFEM) with polygon/quadtree meshes. The PSBFEM is implemented with commercial finite element code Abaqus by the User Element Sub-routine (UEL) feature. The detailed implementation of the framework, defining the UEL element, and solving the stiffness/mass matrix by the eigenvalue decomposition are presented. Several benchmark problems from heat conduction are solved to validate the proposed implementation. Results show that the PSBFEM is reliable and accurate for solving heat conduction problems. Not only can the proposed implementation help engineering practitioners analyze the heat conduction problem using polygonal mesh in Abaqus, but it also provides a reference for developing the UEL to solve other problems using the scaled boundary finite element method.

## 1. Introduction

The analysis of heat conduction has crucial importance in many fields, such as aerospace engineering, mechanical, electrical, and civil engineering. Over the last decades, analytical and numerical methods have been utilized to investigate heat conduction problems. The analytical approaches are limited by the complex geometry and boundary conditions. Numerical methods have become a powerful tool to simulate heat conduction problems.

Due to the equation of transient heat conduction exist space and time variables, we need to discretize the space and the time by the numerical methods. The domain of complex geometry is partitioned into a finite number of non-overlapping subdomains of simplex shapes by introducing the concept of discretization [1]. The most widely used computational techniques in discretized space are the finite element method (FEM) [2]. These methods must discretize the whole computational domain. It is often unavoidable that the generation of highly distorted meshes for complex geometries, which leads to difficulty in solving many complex problems.

In order to overcome the weakness, more alternatives methods are proposed, such as the meshless method [3–5], the smoothed finite element method [6,7], and the polygonal finite element method [7,8]. The scaled boundary finite element method (SBFEM) [9] is a semi-analytical method that attempts to fuse the advantages and characteristics of FEM and the boundary element method (BEM) into one new approach. The SBFEM has been applied to many fields, such as wave propagation [10,11], fracture mechanics [12–15], acoustic [16], seepage [17,18], fluid [19,20], and contract [21]. For these problems, the SBFEM demonstrates more efficiency compared with the conventional FEM.

The polygonal scaled boundary finite element method (PSBFEM) is a novel approach integrating the standard SBFEM and the polygonal mesh technique. This method is flexible in meshing complex geometries, and the use of polygons to discretize the computational domain naturally complements the SBFEM. Polygonal element with more than four sides involves more nodes in their interpolation compared with a conventional FEM. Simultaneously, they generally exhibit superior solution accuracy [1]. Moreover, it is more flexible in the discretization of complex geometry. Therefore, these advantages have further motivated polygonal elements as an alternative to conventional FEM using triangles or quadrilaterals.

Recently, the quadtree mesh generation technique [22][23] is applied to discrete domains. Mesh generation and adaptive refinement of quadtree meshes are straightforward [22]. However, the quadtree discretization will lead to additional nodes, usually called 'handing nodes'. Due to hanging nodes between two adjacent elements of different sizes, it is problematic that quadtree meshes are directly used to simulate within the finite element method's framework. However, the SBFEM only discretizes in the boundary of geometry. Hence, each cell in a quadtree mesh, regardless of hanging nodes, is treated as a generic polygon. The enables the structure of the quadtree to be exploited for efficient computations. The ability to assume any number of sides also enables the SBFEM to discretize curved boundaries better.

In this work, we present a framework for a heat conduction analysis using the SBFEM with polygon/quadtree meshes. We implement steady-state and transient heat conduction analysis of SBFEM by developing a User Element subroutine (UEL). This paper is organized as follows. The governing equations for the heat conduction problem are described in Section 2. The SBFEM of heat conduction problem concept is stated in Section 3. Section 4 outlines the solution procedures. Section 5 describes the implementation of PSBFEM in the heat conduction problems by the Abaqus UEL subroutine. Then, several benchmarks examples are presented in Section 6. Finally, the main concluding remarks of this paper are presented in Section 7.

## 2. Governing equations for the heat conduction problem

In this section, the governing equations for two-dimensional transient heat conduction problems are considered. The governing equations without heat sources are written as

$$\rho c \dot{T} - \nabla \cdot (k \nabla T) = 0 \quad \text{in } \Omega \tag{1}$$

with the initial conditions

$$T(x, y, t = 0) = T_0(x, y) \quad \text{in } \Omega \tag{2}$$

and the boundary conditions

$$T(x, y) = \bar{T} \quad \text{on } S_1 \tag{3}$$

$$-k \frac{\partial T}{\partial n} = \bar{q}_n = \begin{cases} q_2 & \text{on } S_2 \\ h(T - T_\infty) & \text{on } S_3 \end{cases} \tag{4}$$

where $\rho$ is the density, $c$ is the heat capacity, $T$ is the temperature, $\dot{T}$ is the temperature change rate, $\nabla$ is the gradient operator, $k$ is the heat conductivity, $\Omega$ is the computational domain, $n$ is outward normal vector to $\Omega$. $S$ is the boundary. Moreover, $\bar{T}$ and $q_2$ are the prescribed boundary temperature and heat flux, respectively. $h$ is the convection heat transfer coefficient and $T_\infty$ is the ambient temperature.

By applying the Fourier transform on Equation (1), the governing equation is transformed into the frequency domain as

$$k\Delta \tilde{T} = i\omega\rho c \tilde{T} \tag{5}$$

where $\tilde{T}$ is the Fourier transform of $T$. $\omega$ is the frequency. When $\omega=0$, the problems are transformed into the steady-state heat conduction problem.

### 3. The SBFEM for the transient heat conduction problem

As illustrated in Figure 1, the SBFEM presents a local coordinate system $(\xi,\eta)$. The coordinates of a point $(x, y)$ along the radial line and inside the domain can be expressed as follows [24]:

$$x = \xi[N(\eta)]\{x\} \tag{6}$$

$$y = \xi[N(\eta)]\{y\} \tag{7}$$

where $\xi, \eta$ are the scaled boundary coordinates in two-dimensions. $\xi$ is a radial coordinate. $\eta$ is the circumferential coordinate.

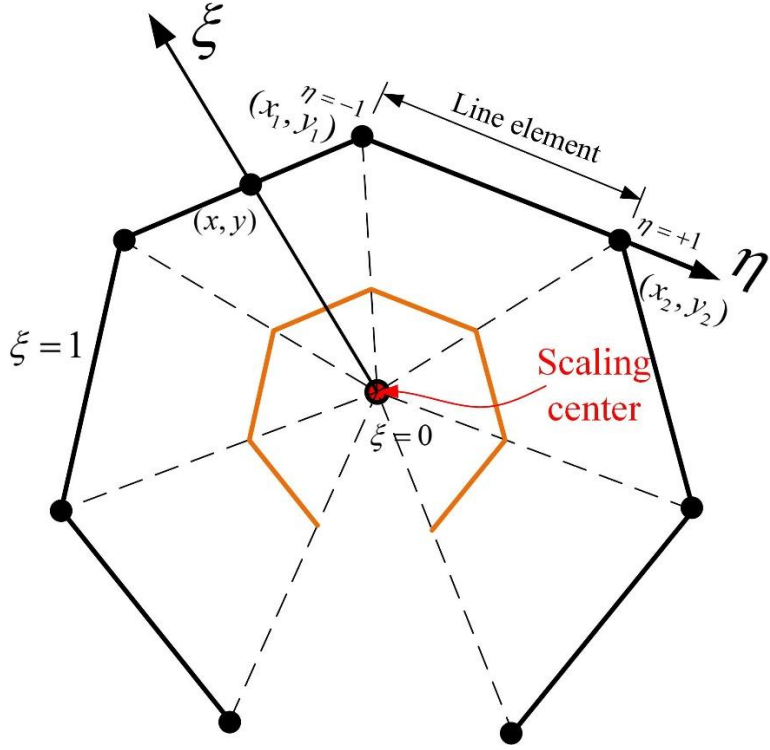

Figure 1. The coordinate system of scaled boundary finite element method.

The differential operator can be transformed from the Cartesian coordinate system $x, y$ into the scaled boundary coordinates system $\xi, \eta$ as follows:

$$\nabla = [b_1]\frac{\partial}{\partial \xi} + \frac{1}{\xi}[b_2]\frac{\partial}{\partial \eta} \tag{8}$$

with

$$[b_1] = \frac{1}{|J_b|}\begin{bmatrix} y_{b,\eta} & 0 \\ 0 & -x_{b,\eta} \\ -x_{b,\eta} & y_{b,\eta} \end{bmatrix} \tag{9}$$

$$[b_2] = \frac{1}{|J_b|}\begin{bmatrix} -y_b & 0 \\ 0 & x_b \\ x_b & -y_b \end{bmatrix} \tag{10}$$

where the Jacobian matrix at the boundary ($\xi = 1$) can be expressed as:

$$[J_b] = \begin{bmatrix} x_b & y_b \\ x_{b,\eta} & y_{b,\eta} \end{bmatrix} = x_b y_{b,\eta} - y_b x_{b,\eta} \tag{11}$$

The temperature $\tilde{T}(\xi,\eta)$ at any point in SBFEM coordinates is written as:

$$\{\tilde{T}(\xi,\eta)\} = [N_u(\eta)]\{\tilde{T}(\xi)\} \tag{12}$$

where $\tilde{T}(\xi)$ is radial temperature functions along a line connecting the scaling center $O$ and a node at the boundary. $[N_u(\eta)]$ is the shape function matrix

$$[N_u(\eta)] = \begin{bmatrix} N_1(\eta) & 0 & N_2(\eta) & 0 & \cdots & 0 & N_M(\eta) & 0 \\ 0 & N_1(\eta) & 0 & N_2(\eta) & \cdots & 0 & 0 & N_M(\eta) \end{bmatrix} \tag{13}$$

By using Equation (8) and (12), the heat flux $q(\xi,\eta)$ can be expressed as

$$q(\xi,\eta) = -k\left([B_1(\eta)]\{\tilde{T}(\xi)\}_{,\xi} + \frac{1}{\xi}[B_2(\eta)]\{\tilde{T}(\xi)\}\right) \tag{14}$$

where,

$$[B_1(\eta)] = \{b_1(\eta)\}[N(\eta)] \tag{15a}$$

$$[B_2(\eta)] = \{b_2(\eta)\}[N(\eta)] \tag{15b}$$

Using the weighted residual method and Green theorem and introducing the boundary conditions, the following SBFEM equation can be written as follows [25]

$$[E_0]\xi^2\{\tilde{T}(\xi)\}_{,\xi\xi} + \left([E_0] + [E_1]^T - [E_1]\right)\xi\{\tilde{T}(\xi)\}_{,\xi} - [E_2]\{\tilde{T}(\xi)\} \\ -[M^0]i\omega\xi^2\{\tilde{T}(\xi)\} = \xi\{F(\xi)\} \tag{16}$$

where,

$$[E_0] = \int_S [B_1(\eta)]^T k[B_1(\eta)]|J_b|\mathrm{d}\eta \tag{17}$$

$$[E_1] = \int_S [B_2(\eta)]^T k[B_1(\eta)]|J_b|\mathrm{d}\eta \tag{18}$$

$$[E_2] = \int_S [B_2(\eta)]^T k[B_2(\eta)]|J_b|\mathrm{d}\eta \tag{19}$$

$$[M_0] = \int_S [N(\eta)]^T \rho c[N(\eta)]|J_b|\mathrm{d}\eta \tag{20}$$

The coefficient matrices $[E_0]$, $[E_1]$, $[E_2]$, and $[M_0]$ that depend only on the geometry and material properties.

## 4. Solution procedure

4.1 Steady-state solution

The SBFEM equation in temperature function $\tilde{T}(\xi)$ is solved for steady-state problems. Setting $\omega=0$ in Equation (16), the SBFEM equations for the temperature field in steady-state heat conduction can be written as

$$[E_0]\xi^2\{\tilde{T}(\xi)\}_{,\xi\xi} + \left([E_0]+[E_1]^T-[E_1]\right)\xi\{\tilde{T}(\xi)\}_{,\xi} - [E_2]\{\tilde{T}(\xi)\} = 0 \quad (21)$$

By introducing the variable,

$$\{X(\xi)\} = \begin{Bmatrix} \{\tilde{T}(\xi)\} \\ \{\tilde{Q}(\xi)\} \end{Bmatrix} \quad (22)$$

The equation can be transformed into first-order ordinary differential equations:

$$\xi\{X(\xi)\}_{,\xi} - [Z_p]\{X(\xi)\} = 0 \quad (23)$$

where the coefficient matrix $[Z_p]$ is a Hamiltonian matrix. The solution for a bounded domain is obtained using the positive eigenvalues of $[Z_p]$. Hence, $[Z_p]$ can be expressed as

$$[Z_p] = \begin{bmatrix} -[E_0]^{-1}[E_1]^T & [E_0]^{-1} \\ [E_2]-[E_1][E_0]^{-1}[E_1]^T & [E_1][E_0]^{-1} \end{bmatrix} \quad (24)$$

The solution of Equation (23) can be obtained by computing the eigenvalue and eigenvector of the matrix $[Z_p]$ yields

$$[Z_p]\begin{bmatrix} [\psi_{11}] & [\psi_{12}] \\ [\psi_{21}] & [\psi_{22}] \end{bmatrix} = \begin{bmatrix} [\psi_{11}] & [\psi_{12}] \\ [\psi_{21}] & [\psi_{22}] \end{bmatrix}\begin{bmatrix} [\lambda_n] & \\ & [\lambda_p] \end{bmatrix} \quad (25)$$

where the real components of eigenvalues $\lambda_n$ and $\lambda_p$ are negative and positive, respectively.

The general solution of equation (23) can be obtained as follows:

$$\{\tilde{T}(\xi)\} = [\psi_{11}]\xi^{-[\lambda_1]}\{c_1\} + [\psi_{12}]\xi^{-[\lambda_2]}\{c_2\} \tag{26}$$

$$\{\tilde{Q}(\xi)\} = [\psi_{21}]\xi^{-[\lambda_1]}\{c_1\} + [\psi_{22}]\xi^{-[\lambda_2]}\{c_2\} \tag{27}$$

where $\{c_1\}$ and $\{c_2\}$ are the integration constants. To obtain a finite solution at the scaling center ($\xi = 0$), $\{c_2\}$ must be equal to zero, the solution in bounded domain can be written as:

$$\{\tilde{T}(\xi)\} = [\psi_{11}]\xi^{-[\lambda_1]}\{c_1\} \tag{28}$$

$$\{\tilde{Q}(\xi)\} = [\psi_{21}]\xi^{-[\lambda_1]}\{c_1\} \tag{29}$$

The relationship between $\{\tilde{T}(\xi)\}$ and $\{\tilde{Q}(\xi)\}$ is expressed as

$$\left[K^{st}\right]\{\tilde{T}(\xi)\} = \{\tilde{Q}(\xi)\} \tag{30}$$

where the steady-state stiffness matrix can be expressed as

$$\left[K^{st}\right] = [\psi_{21}][\psi_{11}]^{-1} \tag{31}$$

4.2 Mass matrix transient solution

To determine the mass matrix $[M]$ of SBFEM, introducing the dynamic-stiffness matrix $[K(\xi,\omega)]$ at $\xi$:

$$[K(\xi,\omega)]\{\tilde{T}(\xi)\} = \{\tilde{Q}(\xi)\} \tag{32}$$

For the bounded domain, the dynamic-stiffness matrix $[K(\omega)]$ on the boundary $\xi=1$ formulated in the frequency domain is written as

$$\left([K(\omega)]-[E_1]\right)[E_0]^{-1}\left([K(\omega)]-[E_1]^T\right) - [E_2] + 2\omega[K(\omega)]_{,\omega} - i\omega[M_0] = 0 \tag{33}$$

To obtain the mass matrix $[M]$ of the bound domain, the low-frequency case in the SBFEM is address, which the dynamic-stiffness $[K(\omega)]$ of a bounded domain is assume as

$$[K(\omega)] = \left[K^{st}\right] + i\omega[M] \qquad (34)$$

The steady-state stiffness matrix $[K^{st}]$ is computed using Equation (31). Substituting Equation (34) into Equation (33) leads to a constant term independent of $i\omega$, a term in $i\omega$ and higher order term in $i\omega$, which are neglected. Besides, the constant term vanishes. The coefficient matrix of $i\omega$ can be expressed as

$$\left(\left(-\left[K^{st}\right]+[E_1]\right)[E_0]^{-1}-[I]\right)[M]+[M]\left([E_0]^{-1}\left(-\left[K^{st}\right]+[E_1]^T\right)-[I]\right)+[M_0]=0 \qquad (35)$$

This is a linear equation to solver the mass matrix $[M]$. By using the eigenvalues and eigenvectors of matrix $[Z_P]$, the Equation (35) can be written as

$$\left([I]+[\lambda_b]\right)[m]+[m]\left([I]+[\lambda_b]\right)=[\psi_{11}]^T[M_0][\psi_{11}] \qquad (36)$$

where,

$$[m]=[\psi_{11}]^T[M][\psi_{11}] \qquad (37)$$

After solving matrix $[m]$ in the Equation (36), the mass matrix $[M]$ is obtained by

$$[M]=\left([\psi_{11}]^{-1}\right)^T[m][\psi_{11}]^{-1} \qquad (38)$$

4.3 Transient solution and time discretization

Substitute Equation (34) into Equation (32) followed by performing inverse Fourier transform, the nodal temperature relationship of a bounded domain is expressed as a standard time domain equation using steady-state stiffness and mass matrices as

$$\left[K^{st}\right]\{T(t)\}+[M]\{\dot{T}(t)\}=\{Q(t)\} \qquad (39)$$

where nodal temperature $T(t)$ is the continuous derivative of time, In general, it is not easy to solve the function solution in the time domain. In this paper, the backward difference method is adopted to solve Equation (39). The time-domain is divided into several time units, and the solution of the time node is obtained step by step from the initial conditions, and the node temperature at any time is obtained by interpolation.

At time $[t, t+\Delta t]$, the temperature change rate $\{\dot{T}(t)\}$ can be expressed as

$$\{\dot{T}(t)\}=\frac{[\Delta T]}{\Delta t}=\frac{\{T(t)\}^{t+\Delta t}-\{T(t)\}^{t}}{\Delta t} \tag{40}$$

Substitute Equation (40) into Equation (39), the equation at time step $t+\Delta t$ can be obtain as follows

$$([K^{st}]^{t+\Delta t}+\frac{[M]^{t+\Delta t}}{\Delta t})\{T(t)\}^{t+\Delta t}=\{Q(t)\}^{t+\Delta t}+\frac{[M]^{t+\Delta t}}{\Delta t}\{T(t)\}^{t} \tag{41}$$

## 5. Framework implementation of PSBEFM for the heat conduction problems

5.1 Overview framework

This work proposes a framework for the heat conduction analysis using the PSBFEM with polygon/quadtree meshes, as shown in Figure 2. The framework contains three sub-modules: the pre-processing modules, the heat conduction analysis modules, and the post-processing modules. In the pre-processing, we develop a python script to generate a polygonal mesh and quadtree mesh. The UEL subroutine is written by FORTRAN 77. Due to the Abaqus CAE not supporting the UEL element's visualization, we extract the 'vtu' format results from the 'odb' file and visualize it in the software Paraview [38].

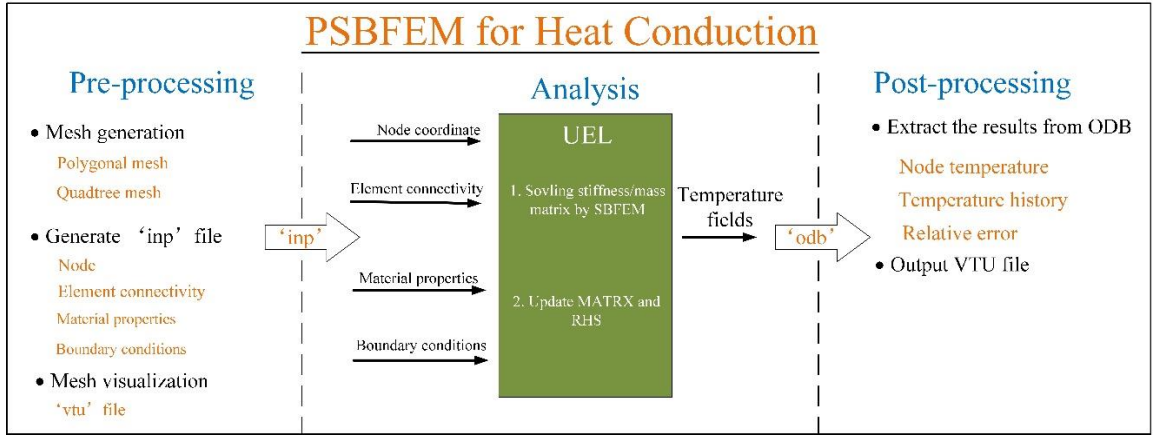

Figure 2. The framework of PSBFEM for heat conduction

## 5.2 The implementation of UEL

The Abaqus/Standard analysis provides a programming interface UEL to define the customized elements. This section presents the essential implementation details of UEL for the SBFEM of heat conduction in Abaqus.

The most critical work of UEL is to update the contribution of the element to the internal force vector RHS and the stiffness matrix AMATRX according to the information from ABAQUS/Standard analysis. Figure 3 shows an overall implementation of the UEL subroutine. Based on the input file's connectivity information, the UEL compute the scaling centers and transform the global coordinate to the local coordinate. The Equation (17)~(20) computes the coefficient matrices $[E_0]$, $[E_1]$, $[E_2]$, and $[M_0]$ which are used to construct the Hamilton matrix $[Z_p]$ using the Equation (24). The two eigenvector matrices ($[\psi_{21}]$, $[\psi_{11}]^{-1}$) are constructed using the eigenvalue decomposition. Further, we can obtain the stiffness matrix $[K]$ and mass matrix $[M]$ of the PSBFEM elements. Finally, the RHS and AMATRX are defined as Equation (30) and Equation (41). In order to decrease the calculation cost, we only solve the stiffness matrix $[K]$, and mass matrix $[M]$ at the first incremental step and store it in the state variable. In the next incremental step, the stiffness/mass matrix is read directly to update RHS and AMATR.

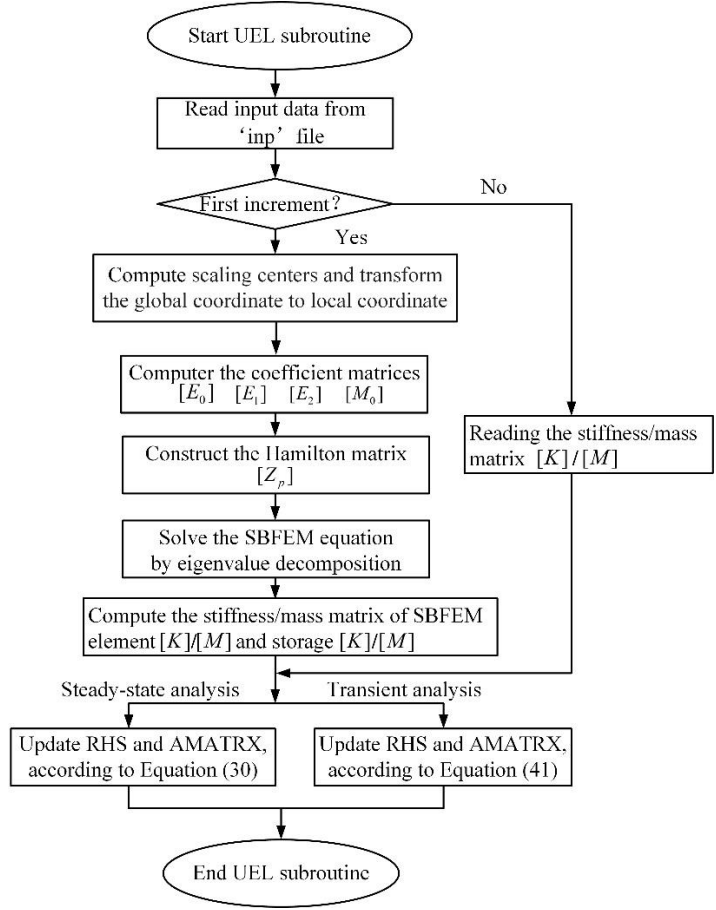

Figure 3. Flowchart of UEL subroutine for SBFEM

To solve the stiffness matrix and mass matrix, we need to employ the eigenvalue decomposition, as shown in Equation (25). At present, many mathematical libraries to perform the eigenvalue decomposition exists. This paper uses the Intel Math Kernel Library (MKL) [36] to decompose the eigenvalue.

5.3 Defining the UEL polygonal elements

The PSBFEM is inherently appropriate for modeling polygons and has other promising capabilities. Because the PSBFEM is discretized only in the boundary, and an element of PSBFEM can assume more complex shapes than a finite element method, the PSBFEM element is much more flexible than the FEM element. With such a significant advantage, the PSBFEM is more applicable to complex geometries than other methods. Figure 4 shows that the supported element type in the PSBFEM. The PSBFEM also provides the complex element type: polygonal and complex quadrilateral elements (quadtree discretization). Hence, it is an effective tool to solve complex elements in numerical analysis.

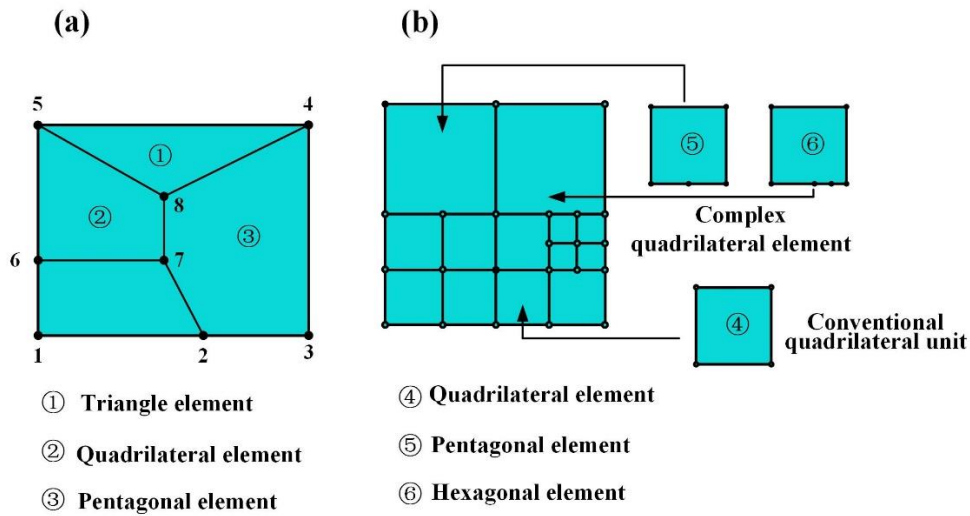

Figure 4. Schematic diagram of the polygon meshes; (a) random polygon meshes; (b) quadtree meshes

The input file of Abaqus contains a complete description of the numerical model, such as node coordinate, element connectivity, the freedom degree, and material properties. This information needs to be defined in the 'inp' file by the user. We show a simple polygon mesh of PSBFEM that demonstrates the definition of elements in the UEL, as shown in Figure 4 (a). In the input file (see Listing 1), the mesh consists of three element types: triangular element (U3), quadrilateral element (U4), and pentagon element (U5). **1 ～ 6** is the line number. Lines **1 ～ 6** are used to define the pentagonal element (U5); Line **1** assigns the element type, the number of nodes, the number of element properties, and the number of freedom degrees per node; Line **2** sets the active degrees of freedom. Lines **3 ～ 4** define the element sets' E5'; Line **6** sets the element properties (conductivity, density, specific heat) of 'E5'.

Listing. 1 The input file of polygon element in the UEL (c.f. Figure 4(a))
```
1 *USER ELEMENT, NODES=5, TYPE=U5, PROPERTIES=2, COORDINATES=2

2 11
3 *ELEMENT, TYPE=U5, ELSET=E5
4 3,2,3,4,8,7
```

```
5 *UEL PROPERTY, ELEST=E5
6 1.0,1.0,1.0
```

## 6. Numerical examples

In this section, we carried out several benchmark problems to demonstrate the convergence and accuracy of the proposed framework for the heat conduction analysis. Moreover, the results of the PSBFEM are compared with the FEM. The FEM analysis uses the commercial finite element software Abaqus. For validation, the relative errors in the temperature are investigated as follows:

$$e_{(\Omega)} = \left\| \mathbf{T} - \mathbf{T}^h \right\|_{(\Omega)} = \frac{\sqrt{\int_{\Omega} \left( \mathbf{T} - \mathbf{T}^h \right)^{\mathrm{T}} \left( \mathbf{T} - \mathbf{T}^h \right) \mathrm{d}\Omega}}{\sqrt{\int_{\Omega} \mathbf{T}^{\mathrm{T}} \mathbf{T} \mathrm{d}\Omega}} \tag{42}$$

where $\mathbf{T}$ is the numerical solution, and $\mathbf{T}^h$ is the analytical or reference solution.

6.1 Steady-state heat conduction analysis

6.1.1 The steady-state heat conduction of square plate

A classical steady-state heat conduction problem is considered to investigate the accuracy and convergence of PSBFEM, as shown in Figure 5. The temperature on the top side is specified as $T = 100 \sin\left(\frac{\pi}{a} x\right)$, and the one at other side is zero. The analytical solution for this problem is given by [26]

$$T(x, y) = \frac{100}{\sinh(0.5\pi)} \sin\left(\frac{\pi}{a} x\right) \cdot \sinh\left(\frac{\pi}{a} y\right) \tag{43}$$

In this example, $a$ is set to 10.0 m. $b$ is set to 5.0m. Due to the symmetry of this problem, we only select the left part of the domain to analyze. The domain is discretized with the quadrilateral and polygonal elements, as shown in Figure 5 (b) and (c). A convergence study is performed by mesh refinement. The meshes are refined successively following the sequence 1m, 0.5m, 0.25m and 0.125m. In this work, the FEM analysis uses the DC2D4 element.

Figure 6 shows the convergence of the relative error in the temperature with mesh refinement. It is observed that all techniques asymptotically converge to an exact solution with an optimal convergence rate (2:1). Moreover, the PSBFEM show slightly accurate results than FEM at the same element size. Figure 7 presents the contours of the temperature of FEM and PSBFEM. It is clearly shown that the results are virtually the same for the FEM and PSBFEM.

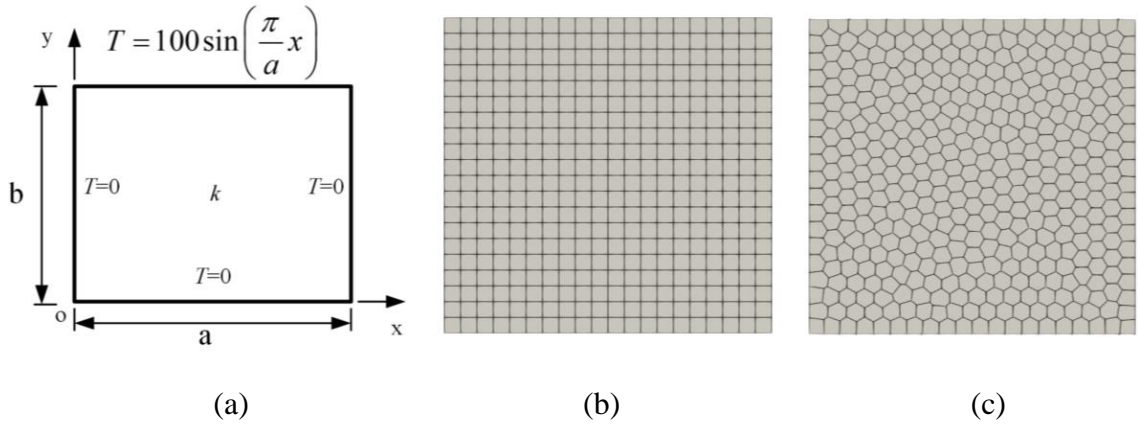

(a)            (b)            (c)

Figure 5. Steady-state heat conduction for the homogeneous material; (a) geometry and boundary conditions; (b) the quadrilateral mesh; (c) the polygonal mesh

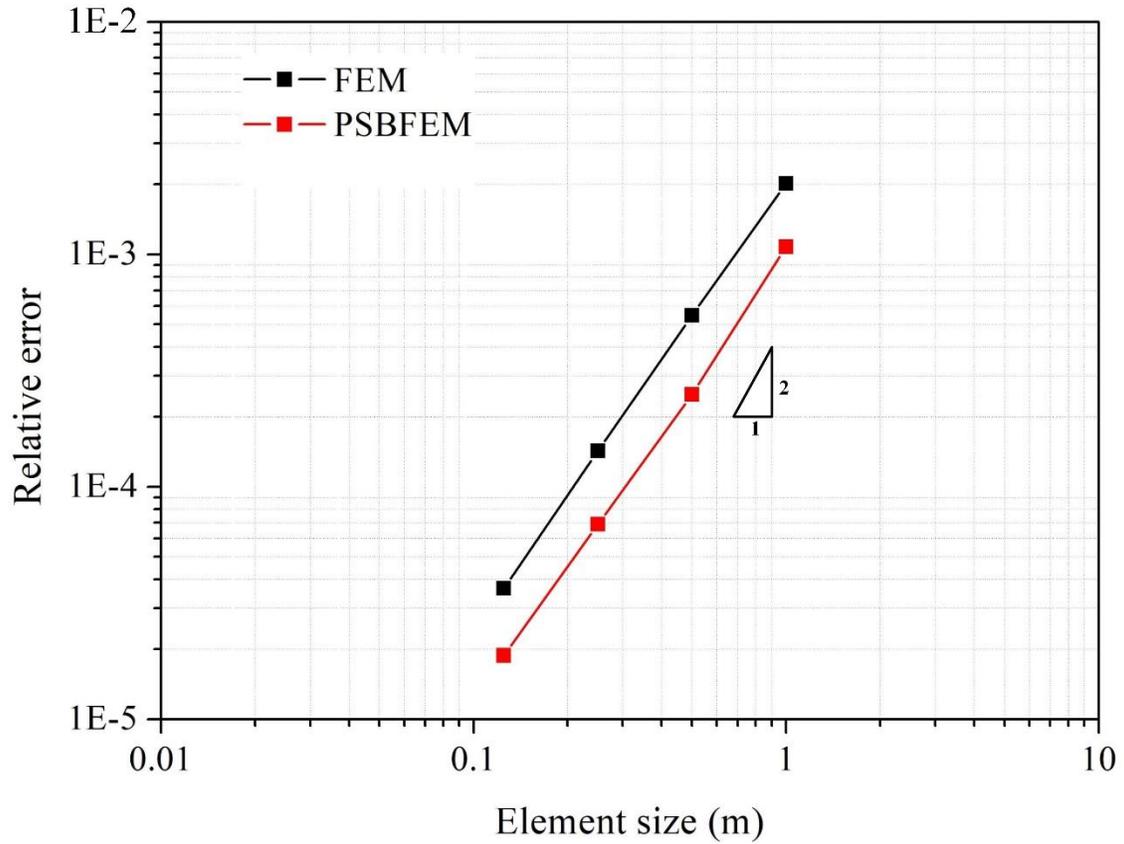

Figure 6. Convergence of the relative error in the temperature

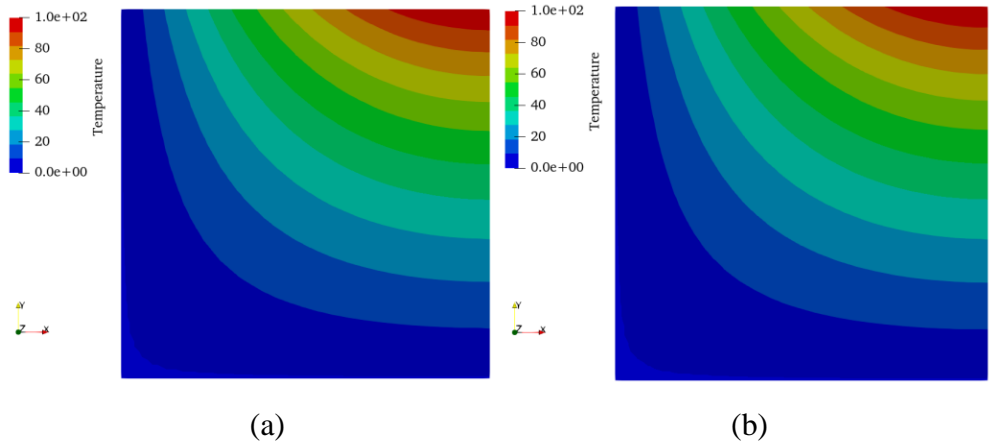

(a)  (b)
Figure 7. Contour plots of temperature; (a) FEM; (b) PSBFEM

6.1.2 The steady-state heat conduction of square plate with multiple holes

In this example, to show the flexibility of the PSBFEM and the quadtree algorithm in the handling mesh, we consider a unit square panel ($L=1.0$ m) with fifteen randomly distributed holes of different sizes, shown in Figure 8. The temperature at the top side is specified as: $T = 150\sin(\frac{\pi}{2}x)$, and the temperature of the bottom side is 10 °C.

The quadtree mesh is generated by setting the same number of mesh seeds on every hole, as shown in Figure 9. It is clearly present that the mesh transition between the holes of different sizes is effectively handled. The square body with multiple holes is modeled with 1632 quadtree elements, and the total nodes are 2478. Moreover, this problem is also analyzed with a similar number of nodes (2532 nodes) using the Abaqus DC2D4 element.

Figure 10 shows that the comparison between the PSBFEM quadtree element and the Abaqus DC2D4 element in the temperature at the right edge. Simple regression is investigated to compare the results of FEM and PSBFEM. Results show that the coefficient of determination $R^2$ is 0.9997, and the relative error is 0.67%. In addition, the temperature distributions obtained from the PSBFEM and the FEM are shown in Figure 11. It is noted that the contour plots present a good agreement. Hence, these results demonstrate PSBFEM accuracy and reliability for the quadtree mesh.

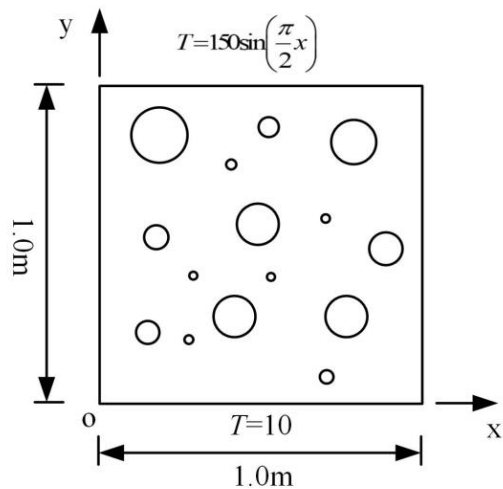

Figure 8. The steady-state heat conduction problem of square panel with multiple holes and boundary conditions

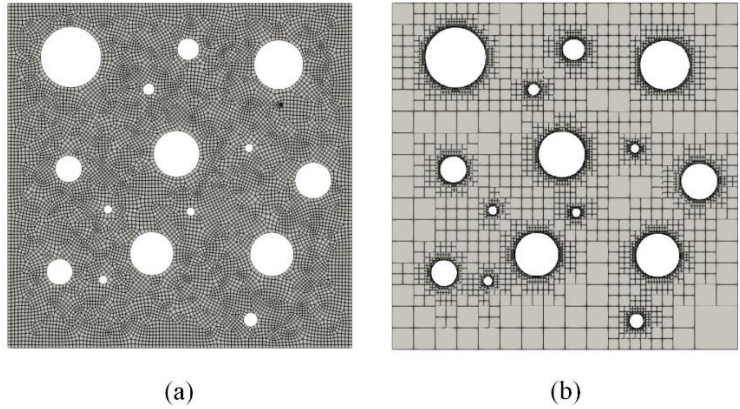

(a)                          (b)

Figure 9. The mesh of square panel with multiple holes; (a) quadrilateral mesh; (b) quadtree mesh

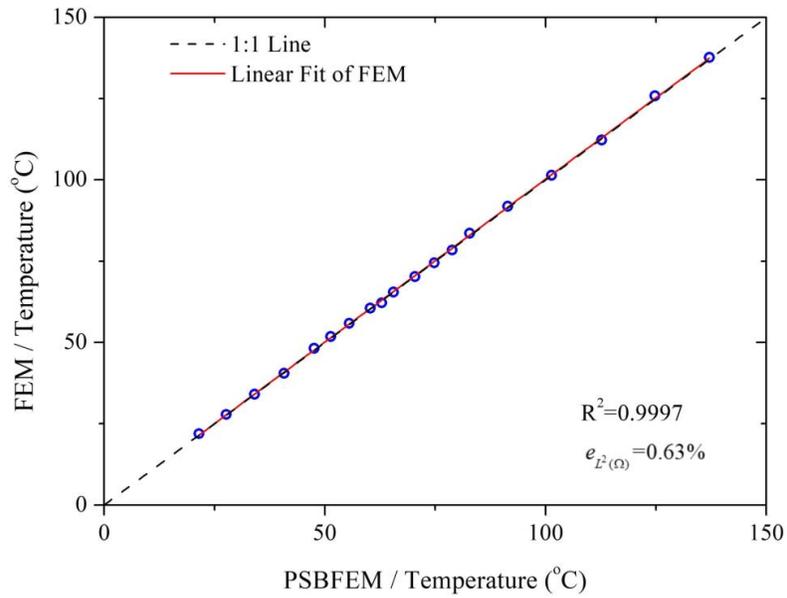

Figure 10. Comparison between quadtree PSBFEM and FEM in the temperature on the right edge

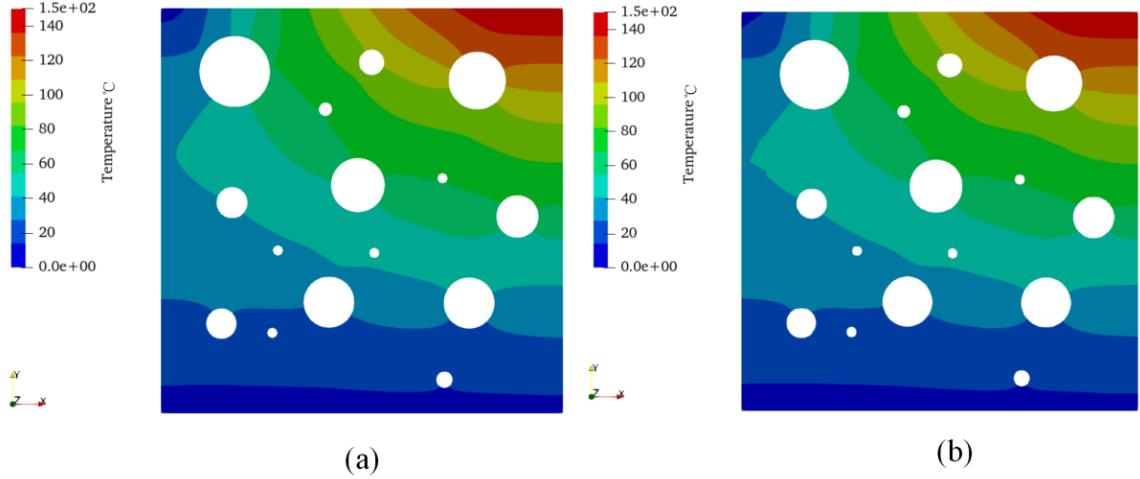

Figure 11. The temperature distribution of multiple holes panel;(a) FEM;(b) PSBFEM

6.2 Transient heat conduction analysis

6.2.1 The transient heat conduction of square plate without heat source

To verify the accuracy of the PSBFEM for the transient heat conduction problem, we consider a square $(\pi \times \pi)$ domain with the analytical solution, as shown in Figure 12. The initial condition is $T_0 = 10\sin(x)\sin(y)$. The Dirichlet boundary conditions are enforced on all sides of the domain, and the Dirichlet boundary conditions are postulated as zero. The material properties are conductivity is 1.0 W/m/°C, $\rho c = 1.0$ J/$\left(\text{m}^3 \cdot {}^\circ\text{C}\right)$. The analytical solution can be expressed as [27]

$$T(x,y,t) = 10e^{-2t}\sin(x)\sin(y) \tag{44}$$

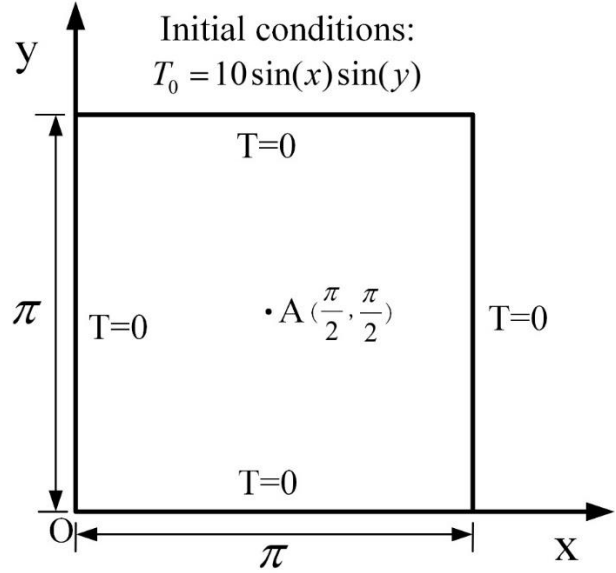

Figure 12. The transient state heat conduction in homogeneous material; (a) geometry and boundary conditions; (b) the quadrilateral mesh; (c) the polygonal mesh

The time step $\Delta t = 0.001s$ and total time $t = 2s$ are used in PSBFEM and FEM. The domain is discretized with the quadrangle and arbitrary polygonal elements. A convergence study is performed by mesh refinement.

The convergence of the relative error in the temperature with mesh refinement is presented in Figure 13. It is observed that the PSBFEM converges to an analytical solution with an optimal convergence rate. The temperature-time history of point A by the FEM and PSBFEM is compared in Figure 14, where the results obtained by the two methods correspond well with the analytical solution. In addition, the history temperature distributions at different times are presented in Figure 15. It is noted that the contour plots present a good agreement. Therefore, these results demonstrate that the PSBFEM is reliable and accurate for solving the transient heat conduction problems.

To evaluate the time consumption of the PSBFEM, we compare the PSBFEM and Abaqus standard elements using the quadrilateral mesh. Since the increment size setting affects the solution time, we use the auto-increment type. The comparison is evaluated with an Intel Core i7-4710MQ CPU running at 2.50GHz and 4.0 GB of RAM. The total CPU time is normalized. Results presented in Figure 16 show that the time consumption of PSBFEM and Abaqus standard elements are comparable.

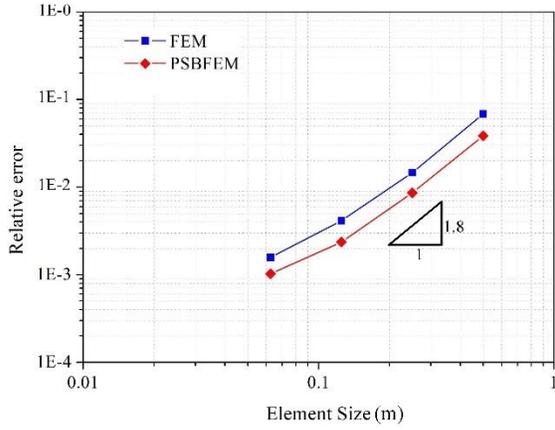 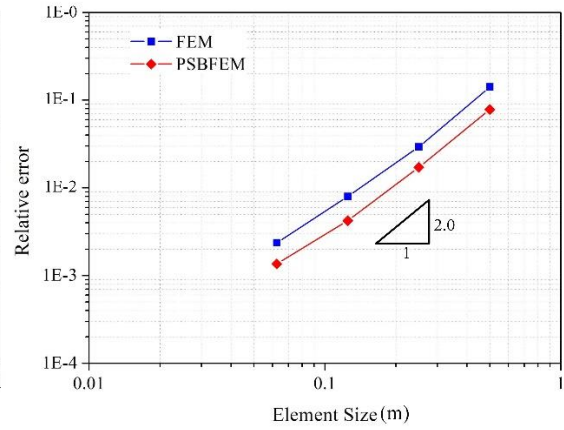

(a) T=0.5s  (b) T=1.0s

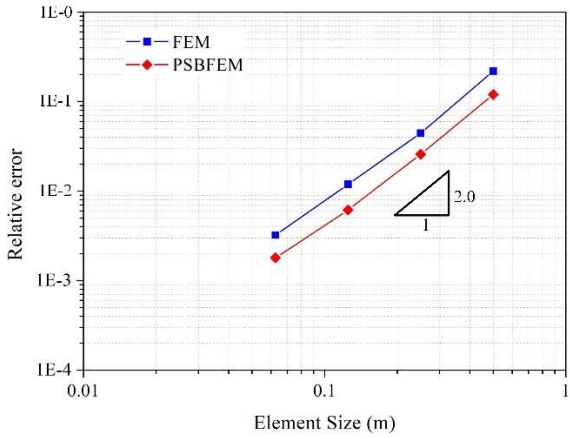 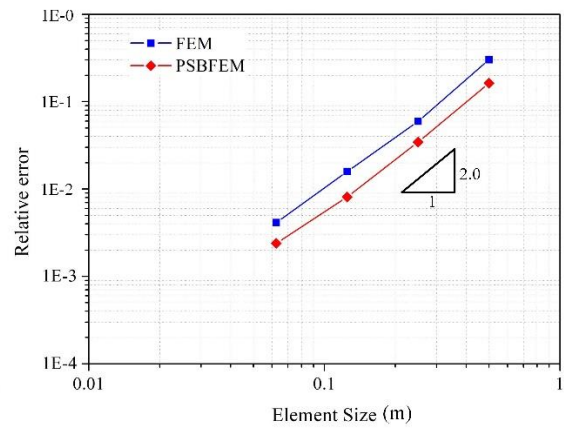

(b) T=1.5s  (d) T=2.0s

Figure 13. Convergence of the relative error in the temperature at the different times

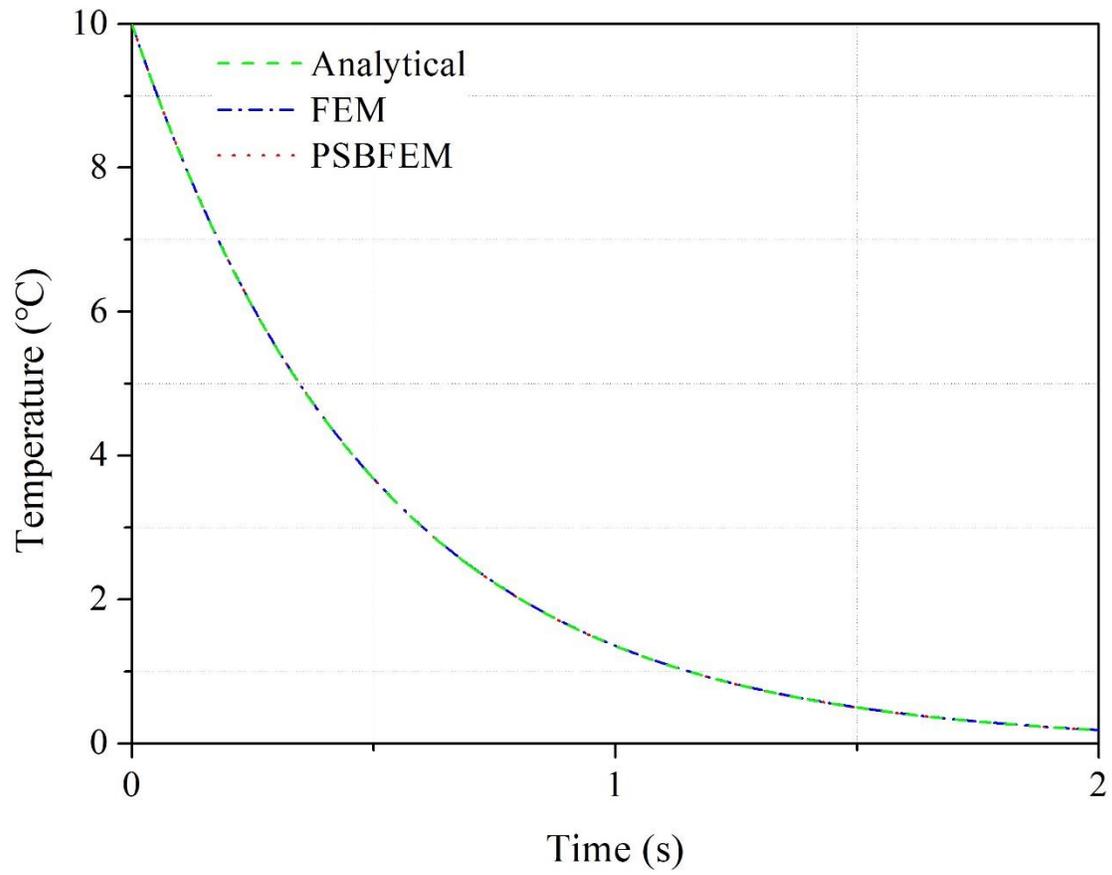

Figure 14. Comparison of temperature time history of point A

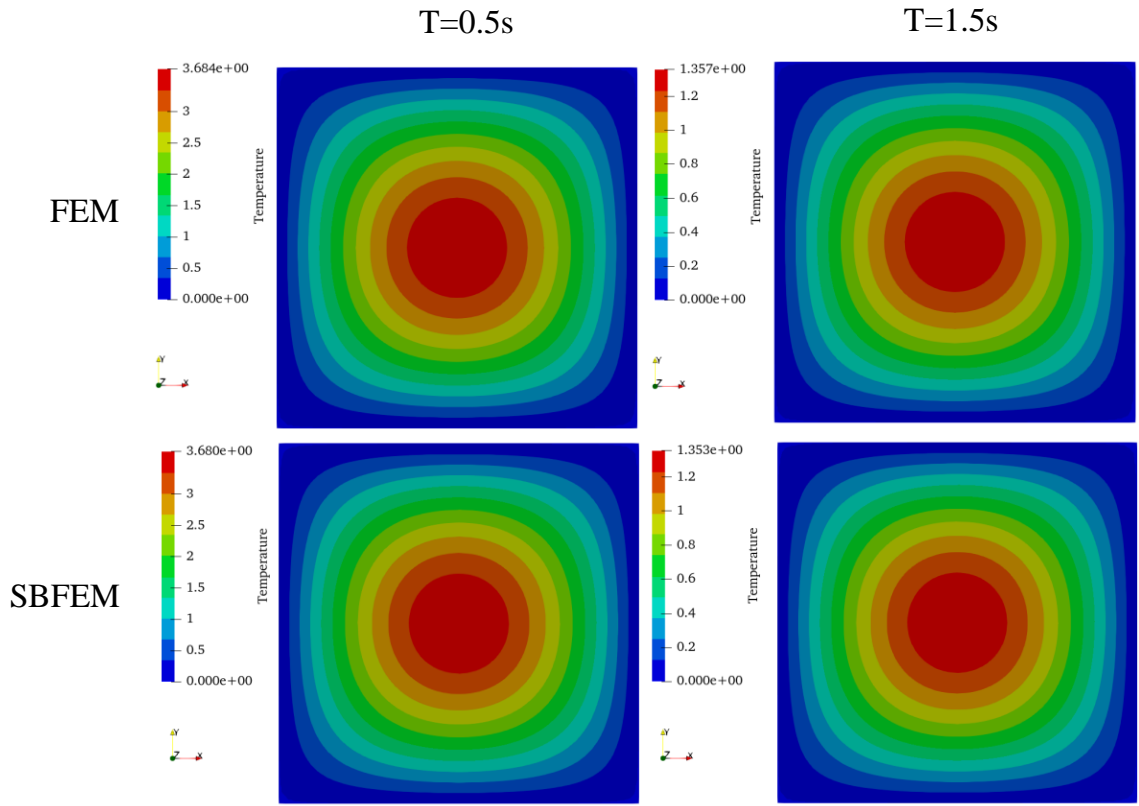

Figure 15. The history temperature distribution at different time using FEM and PSBFEM

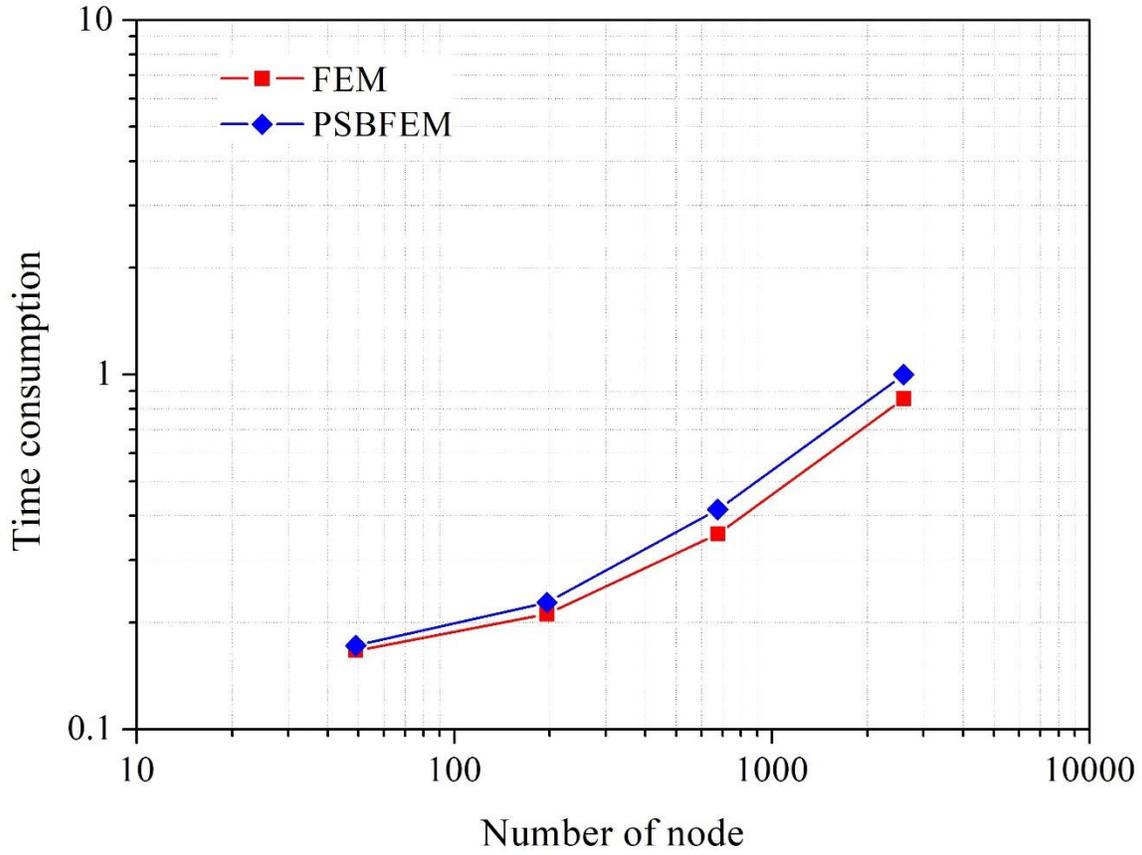

Figure 16. The time consumption comparison of PSBFEM between Abaqus standard elements

6.2.2 The transient heat conduction for the complex geometry

In this example, we consider a square plate ($L=4.0m$) with a rabbit shape cavity, as shown in Figure 17. The material properties are conductivity is 52W/m/°C, specific heat is 434J/kg/°C, density is 7800kg/m$^3$. To compare results, four nodes A, B, C, D are chosen as shown in Figure 17. The temperature at the top side is specified as $300sin(\frac{\pi}{4}x)$, and the one at other side is 50°C. The time step $\Delta t = 20s$ and total time $t = 1000s$ are used in PSBFEM and FEM.

The quadtree mesh is generated by setting the same number of mesh seeds on the cavity, as shown in Figure 18. It is clearly present that the mesh transition between the cavity of different sizes is effectively handled. The square body with a rabbit shape cavity is modelled with 3147 quadtree elements. Moreover, this problem is also analyzed with a similar number of elements using the Abaqus DC2D4 element, and the total elements are 3358.

It is noted that the PSBFEM obtained solutions are in a very good agreement with the Abaqus for the all points, as shown in Figure 19. Table. 2 shows that the relative error of four nodes is less than 1.0%. Moreover, Figure 20 presents that the temperature distribution at different time using FEM and PSBFEM. The temperature distribution are virtually the same for the FEM and PSBFEM. Therefore, the PSBFEM with quadtree meshes shows a good effect for solving complex geometric in the transient heat conduction problem.

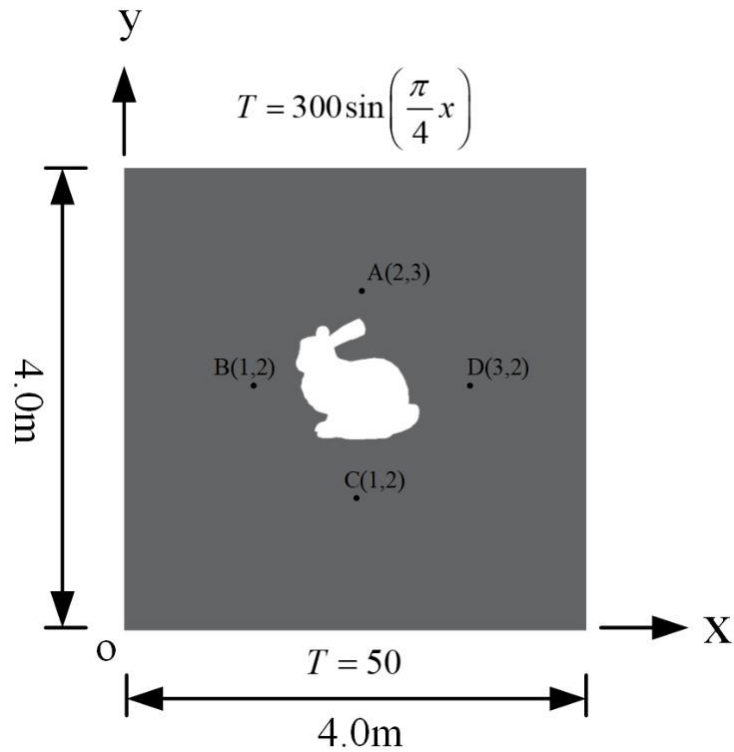

Figure 17. The geometry and boundary conditions of square plate with a rabbit shape cavity

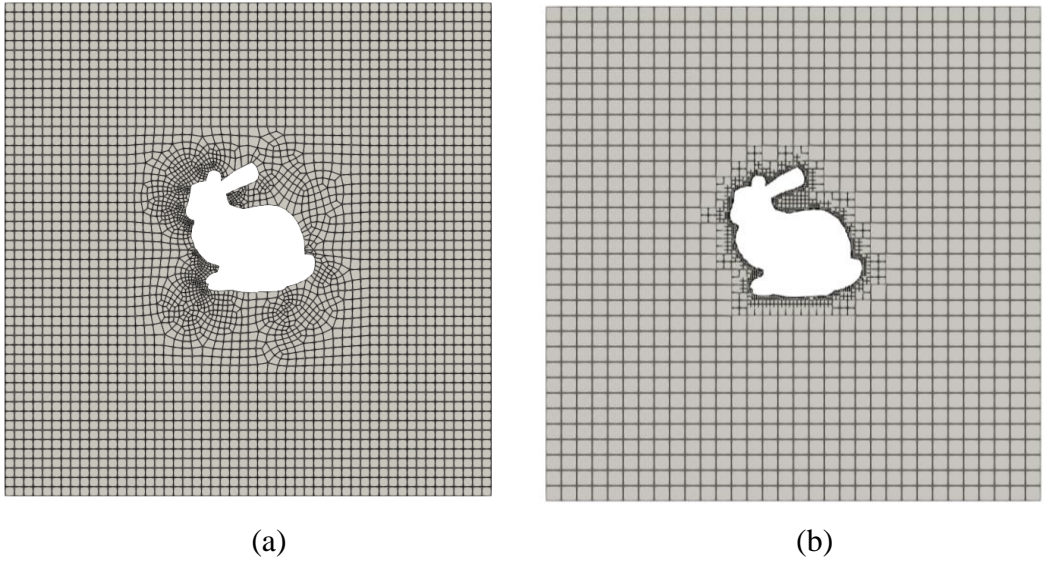

(a)                          (b)

Figure 18. The meshes of square plate with a rabbit shape cavity; (a) Abaqus mesh (b) quadtree mesh

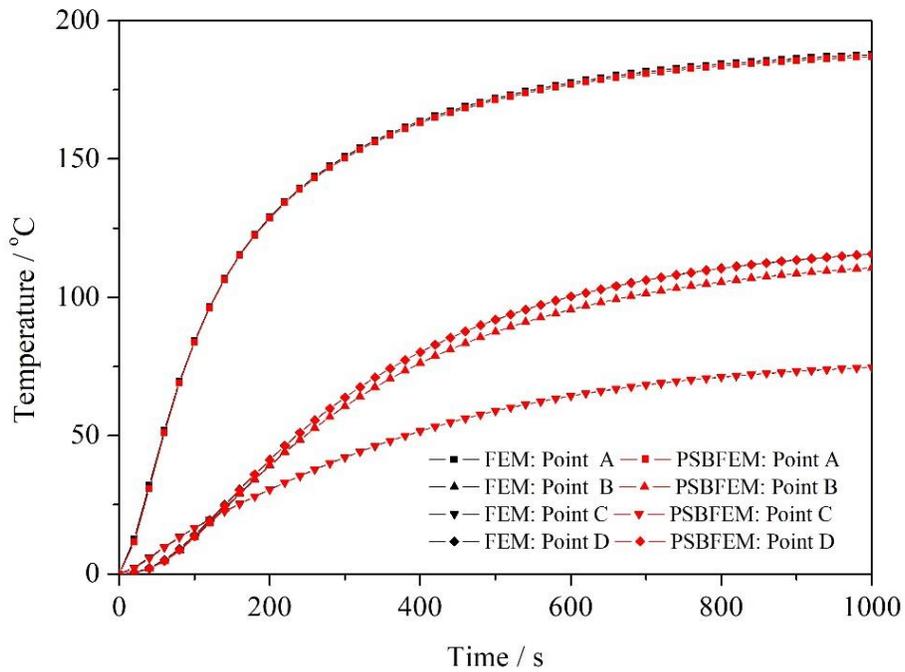

Figure 19. The history temperature of four nodes

Table. 2 The relative error for the FEM and PSBFEM

| Points | A | B | C | D |
|---|---|---|---|---|
| Relative error (%) | 0.42 | 0.15 | 0.14 | 0.15 |

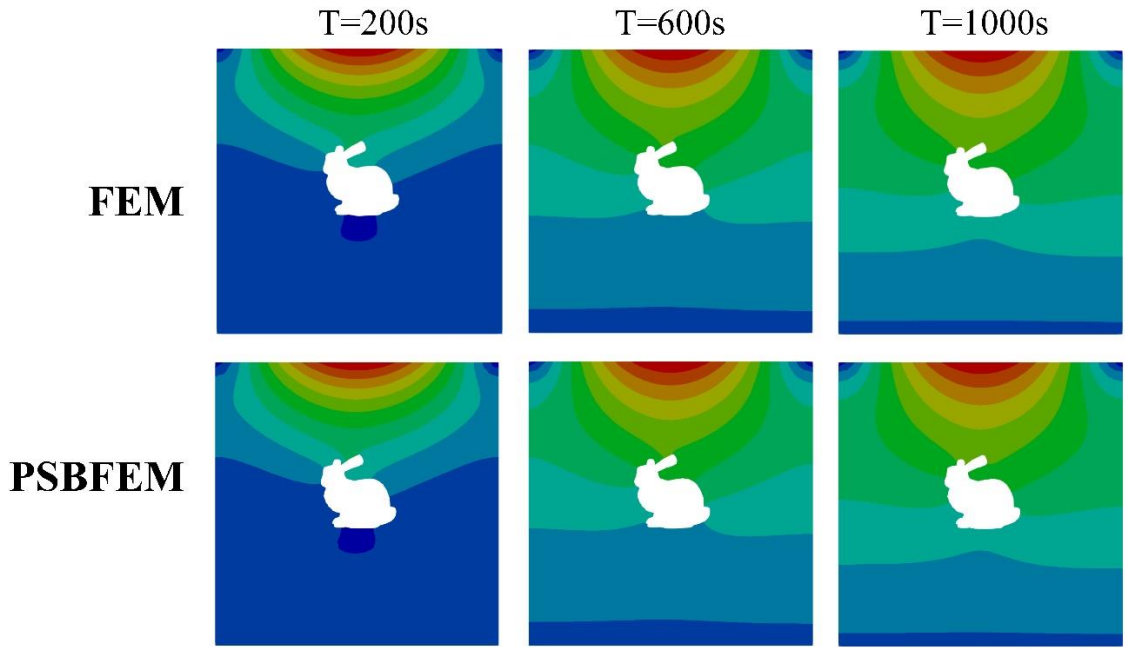

Figure 20. The temperature distribution at different time using FEM and PSBFEM.

## 7. Conclusions

This paper implements a framework for the heat conduction problems using the SBFEM and the polygon/quadtree discretization within the Abaqus/Standard analysis by the UEL subroutine. This work mainly focuses on the detailed implementation of the framework, input data format, and the UEL subroutine, one of the proposed work's critical features.

The implementation of PSBFEM is validated against the FEM by solving a several benchmark problems. The polygon element of PSBFEM has a higher accuracy rate than the standard FEM element. In addition, the results also show that PSBFEM converges to an analytical solution with an optimal convergence rate. The proposed framework is robust accurate for solving the steady-state and transient heat conduction problems. In the future, we will extend the new technique to three-dimensional heat conduction problems.